\documentclass[number,citesort,seceqn,dvips]{arxbj}

\aid{0}
\volume{17}
\issue{1}
\pubyear{2011}
\firstpage{155}
\lastpage{169}
\doi{10.3150/00-BEJ268}

\makeatletter

\newtheorem{theorem}{Theorem}[section]

\newtheorem{lem}[theorem]{Lemma}
\newtheorem{prop}[theorem]{Proposition}
\newproclaim{definition}[theorem]{Definition}
\newremark{remarks}[theorem]{Remark}

\def\eqref#1{(\ref{#1})}

\makeatother

\begin{document}
\begin{frontmatter}

\title{Transportation inequalities: From Poisson to Gibbs measures}
\runtitle{Transportation inequalities}

\begin{aug}
\author[a]{\fnms{Yutao} \snm{Ma}\thanksref{a}\ead[label=e1]{mayt@bnu.edu.cn}},
\author[b]{\fnms{Shi} \snm{Shen}\thanksref{b}\ead[label=e2]{henlen80@yahoo.com.cn}},
\author[c]{\fnms{Xinyu} \snm{Wang}\thanksref{c}\ead[label=e3]{wang\_xin\_yu2000@hotmail.com}}
\and
\author[d]{\fnms{Liming} \snm{Wu}\corref{}\thanksref{d}\ead[label=e4]{Li-Ming.Wu@math.univ-bpclermont.fr}}
\runauthor{Ma, Shen, Wang and Wu }
\address[a]{School of Mathematical Sciences and Lab. Math. Com. Sys.,
Beijing Normal University, 100875 Beijing, China. \printead{e1}}
\address[b]{College of Science, Minzu University of China,
100081 Beijing, China.\\ \printead{e2}}
\address[c]{School of Mathematics, Wuhan University, 430072 Hubei,
China.\\ \printead{e3}}
\address[d]{Laboratoire de Math\'ematiques, CNRS UMR 6620, Universit\'
e Blaise
Pascal, avenue des Landais 63177 Aubi\`ere, France and Institute
of Applied Mathematics, Chinese Academy of Sciences, 100190 Beijing,
China. \printead{e4}}
\end{aug}

\received{\smonth{3} \syear{2009}}
\revised{\smonth{2} \syear{2010}}

%
\begin{abstract}
We establish an optimal transportation inequality for
the Poisson measure on the configuration space. Furthermore, under
the Dobrushin uniqueness condition, we obtain a sharp transportation
inequality for the Gibbs measure on $\mathbb{N}^\Lambda$ or the
continuum Gibbs
measure on the configuration space.
\end{abstract}

%

\end{frontmatter}

\section{Introduction}\label{sec1}

 \textit{Transportation inequality $W_1H$.} Let $\mathcal{X}$ be a
Polish space equipped with the Borel $\sigma$-field $\mathcal{B}$ and $d$
be a lower semi-continuous metric on the product space $\mathcal
{X}\times
\mathcal{X}$ (which does not necessarily generate the topology of
$\mathcal{X}$).
Let $\mathcal{M}_1(\mathcal{X})$ be the space of all probability
measures on $\mathcal{X}$.
Given $p\ge1$ and two probability measures $\mu$ and $\nu$ on
$\mathcal{X}$, we define the quantity
\[
W_{p, d}(\mu, \nu)=\inf\biggl(\int \hspace{-2pt}\int  {d}(x, y)^p\,\mathrm{d}\pi(x, y)\biggr)^{1/p},
\]
where the infimum is taken over all probability measures $\pi$ on
the product space $\mathcal{X}\times\mathcal{X}$ with marginal
distributions $\mu$
and $\nu$ (say, coupling of $(\mu, \nu)$). This infimum is finite
provided that $\mu$ and $\nu$ belong to $\mathcal{M}_1^p(\mathcal
{X},d):=\{\nu\in
\mathcal{M}_1(\mathcal{X});
\int d^p(x,x_0)\,\mathrm{d}\nu<+\infty\}$, where $x_0$ is some fixed point of
$\mathcal{X}$.
This quantity is commonly referred
to as the \textit{$L^p$-Wasserstein distance} between $\mu$ and~$\nu
.$ When
$d$ is the trivial metric $d(x, y)= 1_{x\neq y},  2W_{1,d}(\mu,
\nu)=\|\mu-\nu\|_{\mathrm{TV}},$ the total variation of $\mu-\nu.$

The Kullback information (or relative entropy) of $\nu$ with respect
to $\mu$ is defined as
%
\begin{equation}\label{kullback} H(\nu/\mu)=
\cases{
\displaystyle\int\log\dfrac{\mathrm{d}\nu}{\mathrm{d}\mu}\,\mathrm{d}\nu &\quad\mbox{if}  $\nu\ll \mu$, \cr
+\infty&\quad\mbox{otherwise}.
}
\end{equation}
Let $\alpha$ be a non-decreasing left-continuous function on
$\mathbb{R}^+=[0,+\infty)$ which vanishes at $0$. If, moreover,
$\alpha$ is
convex, we write $\alpha\in\mathcal{C}$. We say that the probability measure
$\mu$ satisfies the \textit{transportation inequality $\alpha$-$W_1H$
with deviation function $\alpha$} on $(\mathcal{X}, d)$ if
%
\begin{equation}\label{W1H*} \alpha(W_{1, d}(\mu, \nu)
)\le H(\nu
/\mu)
\qquad  \forall\nu\in\mathcal{M}_1(\mathcal{X}).
\end{equation}
This transportation
inequality $W_1H$ was introduced and studied by Marton
\cite{Mar96} in relation with measure concentration, for quadratic
deviation function $\alpha$. It was further characterized by
Bobkov and G\"otze \cite{BG99}, Djellout, Guillin and Wu \cite{DGW}, Bolley
and Villani \cite{BV05} and others. The latest development is due to
Gozlan and L\'eonard \cite{GL}, in which the general $\alpha$-$W_1H$
inequality above was introduced in relation to large deviations
and characterized by concentration inequalities, as follows.

\begin{theorem}[(Gozlan and L\'eonard \cite{GL})]\label{Gozleo}
Let $\alpha\in\mathcal{C}$ and $\mu\in\mathcal{M}_1^1(\mathcal
{X},d)$. The
following statements are then equivalent:
\begin{enumerate}[(b$'$)]
\item[(a)] the transportation
inequality $\alpha$-$W_1H$ (\ref{W1H*}) holds;

\item[(b)] for all $\lambda\geq0$
and all $F\in b\mathcal{B}$, $\|F\|_{{\rm Lip}(d)}:=\sup_{x\ne
y}\frac{|F(x)-F(y)|}{d(x,y)}\le1$,
\[
\log\int_\mathcal{X} \exp\bigl(\lambda
[F-\mu(F)]\bigr)\mu(\mathrm{d}x)\leq\alpha^\ast(\lambda),
\]
where $\mu(F):=\int_\mathcal{X}
F\,\mathrm{d}\mu$ and $\alpha^\ast(\lambda):=\sup_{r\ge0}(\lambda r -
\alpha(r))$ is the semi-Legendre transformation of $\alpha$;

\item[(b$'$)] for all $\lambda\geq0$ and all $F,G\in C_b(\mathcal
{X})$ (the
space of all bounded and
continuous functions on $\mathcal{X}$) such that $F(x)-G(y)\le d(x,y)$
for all $x,y\in\mathcal{X}$,
\[
\log\int_\mathcal{X} \mathrm{e}^{\lambda
F}\mu(\mathrm{d}x)\leq\lambda\mu(G)+ \alpha^\ast(\lambda);
\]
\item[(c)] for any measurable function $F$ such that $\|F\|_{{\rm
Lip}(d)}\le1$, the following concentration inequality holds true:
for all $n\geq1, r\geq0$,
%
\begin{equation}\label{gozleo-1}
\mathbb{P}\Biggl(\frac{1}{n}\sum
_1^nF(\xi_k)\geq\mu(F)+r\Biggr)\leq
\mathrm{e}^{-n\alpha(r)},
\end{equation}
where $(\xi_n)_{n\geq1}$ is a
sequence of i.i.d.~$\mathcal{X}$-valued random variables with common
law $\mu$.
\end{enumerate}
\end{theorem}

The estimate on the Laplace transform in (b) and the concentration
inequality in (\ref{gozleo-1}) are the main motivations for the
transportation inequality ($\alpha$-$W_1H$).

\textit{Objective and organization}. The objective
of this paper is to prove the transportation inequality
$(\alpha$-$W_1H)$ for:
\begin{enumerate}[(1)]
\item[(1)] (the free case) the Poisson measure $P^0$ on the configuration space
consisting of Radon point measures $\omega=\sum_i \delta_{x_i},
x_i\in E$ with some $\sigma$-finite intensity measure $m$ on $E$,
where $E$ is some fixed locally compact space;

\item[(2)] (the interaction case) the continuum Gibbs measure over a
compact subset $E$ of $\mathbb{R}^d$,
\[
P^{\phi}(\mathrm{d}\omega)=\frac{\mathrm{e}^{-(1/2)\sum_{x_i, x_j\in\operatorname{supp}\omega, i\ne j}\phi(x_i-x_j)
-\sum_{k, x_i\in\operatorname{supp}(\omega)}\phi(x_i-y_k)}}{Z}P^0(\mathrm{d}\omega),
\]
where $\phi\dvtx \mathbb{R}^d\to[0,+\infty]$ is some pair-interaction
non-negative even function (see Section~\ref{sec4} for notation) and $P^0$ is the
Poisson measure with intensity $z\,\mathrm{d}x$ on $E$.
\end{enumerate}

For Poisson measures on $\mathbb{N}$, Liu \cite{Liu09} obtained the optimal
deviation function by means of Theorem \ref{Gozleo}. For
transportation inequalities of Gibbs measures on discrete sites, see
 \cite{Marton04} and  \cite{Wu06}.

For an illustration of our main result (Theorem \ref{w1hgibbs}) on
the continuum Gibbs measure $P^\phi$, let $E:=[-N,N]^d$ ($1\le N\in
\mathbb{N}$) and $f\dvtx [-N,N]^d\to\mathbb{R}$ be measurable and
periodic with
period $1$ at each variable so that $|f|\le M$. Consider the
empirical mean per volume $F(\omega):=\omega(f)/(2N)^d$ of $f$.
Under Dobrushin's uniqueness condition $D:=z
\int_{\mathbb{R}^d}(1-\mathrm{e}^{-\phi(y)})\, \mathrm{d}y<1$, we have (see Remark
\ref{rem41} for proof)
%
\begin{equation}\label{a1}
P^\phi\bigl(F>P^\phi(F) +r \bigr)\le\exp\biggl(-
\frac{(2N)^d(1-D)r}{2M} \log\biggl( 1+\frac{(1-D)r}{zM} \biggr)
\biggr),\qquad  r>0,
\end{equation}
an explicit Poissonian concentration inequality which is sharp when
$\phi=0$.

The paper is organized as follows. In the next section, we prove
$(\alpha$--$W_1H)$ for the Poisson measure on the configuration
space with respect to two metrics: in both cases,
we obtain optimal deviation functions. Our main tool is
Gozlan and Leonard's Theorem \ref{Gozleo} and a known concentration
inequality in \cite{Wu00}. Section~\ref{sec3}, as a prelude to the study of
the continuum Gibbs measure $P^\phi$ on the configuration space, is
devoted to the study of a Gibbs measure on $\mathbb{N}^{\Lambda}$. Our
method is a combination of a lemma on $W_1H$ for mixed measure,
Dobrushin's uniqueness condition and the McDiarmid--Rio martingale
method for dependent tensorization of the $W_1H$-inequality. Finally, in
the last section, by approximation, we obtain a sharp
$(\alpha$--$W_1H)$ inequality for the continuum Gibbs measure
$P^\phi$ under Dobrushin's uniqueness condition $D=z
\int_{\mathbb{R}^d}(1-\mathrm{e}^{-\phi(y)})\,\mathrm{d}y<1$. The latter is a sharp
sufficient condition, both for the analyticity of the pressure
functional and for the spectral gap; see \cite{Wu04}.

\section{Poisson point processes}\label{sec2}

 \textit{Poisson space.} Let $E$ be a metric complete locally
compact space with the Borel field $\mathcal{B}_E$ and $m$ a
$\sigma$-finite positive Radon measure on $E$. The Poisson space
$(\Omega,\mathcal{F},P^0)$ is given by:
\begin{enumerate}[(1)]
\item[(1)] $\Omega:=\{\omega=\sum_i\delta_{x_i} \mbox{(Radon
measure); } x_i\in
E\}$ (the so-called configuration space over $E$);
\item[(2)] $\mathcal{F}=\sigma(\omega\rightarrow\omega(B)|B\in
\mathcal{B}_E)$;
\item[(3)] $\forall B\in\mathcal{B}_E, \forall
k\in\mathbb{N}\mbox{: } P^0(\omega
\dvtx \omega(B)=k)=\mathrm{e}^{-m(B)}\frac{m(B)^k}{k!}$;
\item[(4)] $\forall B_1,\dots,B_n\in\mathcal{B}_E$ disjoint,
$\omega(B_1),\ldots,\omega(B_n)$ are $P^0$-independent,
\end{enumerate}
where $\delta_x$ denotes the Dirac measure at $x$. Under $P^0$,
$\omega$ is exactly the Poisson point process on $E$ with intensity
measure $m(\mathrm{d}x)$. On $\Omega$, we consider the vague convergence
topology, that is, the coarsest topology such that $\omega\to
\omega(f)$ is continuous, where $f$ runs over the space $C_0(E)$ of
all continuous functions with compact support on $E$. Equipped with
this topology, $\Omega$ is a Polish space and this topology is the
weak convergence topology (of measures) if $E$ is compact.

\begin{definition}
Letting $\varphi$ be a positive measurable function on $E$, we define a
metric $d_\varphi(\cdot, \cdot)$ (which may be infinite) on the Poisson
space $(\Omega, \mathcal{F}, P^0)$ by
\begin{eqnarray*}
d_\varphi(\omega, \omega')=\int_E\varphi \,\mathrm{d}|\omega-\omega'|,
\end{eqnarray*}
where $|\nu|:=\nu^++\nu^-$ for a signed measure $\nu$ ($\nu^\pm$ are,
respectively, the positive and negative parts of $\nu$ in the
Hahn--Jordan decomposition).
\end{definition}

\begin{lem}
If $\varphi$ is continuous, then the metric $d_\varphi$ is lower
semi-continuous on $\Omega$.
\end{lem}

\begin{pf} Indeed, for any $\omega,\omega'\in\Omega$,
\[
d_\varphi(\omega, \omega')=\sup_{f} |\omega(f)-\omega'(f)|,
\]
where the supremum is taken over all bounded $\mathcal{B}_E$-measurable
functions $f$ with compact support such that $|f|\le\varphi$. Now,
as $\varphi$ is continuous, we can approximate such $f$ by $f_n\in
C_0(E)$ in $L^1(E, \omega+\omega')$ and $|f_n|\le\varphi$. Then
\[
d_\varphi(\omega, \omega')=\sup_{f\in C_0(E), |f|\le\varphi}
|\omega(f)-\omega'(f)|.
\]
As $(\omega,\omega')\to|\omega(f)-\omega'(f)|$ is continuous on
$\Omega\times\Omega$, $d_\varphi(\omega, \omega')$ is lower
semi-continuous on $\Omega\times\Omega$.
\end{pf}

Assume from now on that $\varphi$ is continuous. Then, for
any $\nu,\mu\in\mathcal{M}_1(\Omega),$ we have the Kantorovitch--Rubinstein
equality \cite{Kel84,Leo,Villani},
\begin{eqnarray*}
W_{1,d_\varphi}(\mu,\nu)
&=&\sup\biggl\{\int F \,\mathrm{d}\nu-\int G\,\mathrm{d}\mu\Big| F,G\in C_b(\Omega),
F(\omega)-G(\omega')\leq d_\varphi(\omega, \omega')\biggr\}
\\
&=&\sup\biggl\{\int G \,\mathrm{d}(\nu-\mu)\dvtx G\in b\mathcal{F},   \|G\|_{{\rm
Lip}(d_\varphi)}\leq1\biggr\}.
\end{eqnarray*}
Here, $b\mathcal{F}$ is the space of all real, bounded and $\mathcal
{F}$-measurable
functions.

\textit{The difference operator $D$.} We denote by
$L^0(\Omega,P^0)$ the space of all $P^0$-equivalent classes of real
measurable functions w.r.t.~the completion of $\mathcal{F}$ by $P^0$. Hence,
the difference operator $D\dvtx L^0(\Omega,P^0)\rightarrow L^0(E\times
\Omega, m\otimes P^0)$ given by
\[
F\rightarrow
D_xF(\omega):=F(\omega+\delta_x)-F(\omega)
\]
is well defined (see \cite{Wu00}) and plays a crucial role in
the Malliavin calculus on the Poisson space.

\begin{lem}\label{lem22} Given a measurable function
$F\dvtx \Omega\rightarrow\mathbb{R}$, $\| F\|_{{\rm
Lip}(d_\varphi)}\leq1$ if and only if $|D_xF(\omega)|\leq\varphi(x)$
for all $\omega\in\Omega$ and $x\in E$.
\end{lem}

\begin{pf} If $\|F\|_{\mathrm{Lip}(d_\varphi)}\le1$,
since
\[
|D_xF(\omega)|=|F(\omega+\delta_x)-F(\omega)|\leq
d_\varphi(\omega+\delta_x, \omega)=\int_E \varphi
\,\mathrm{d}|(\omega+\delta_x)-\omega|=\varphi(x),
\]
the necessity is true. We now prove the sufficiency. For any $\omega,
\omega'\in\Omega,$ we
write $\omega=\sum_{k=1}^i\delta_{x_k}+\omega\wedge\omega'$ and
$\omega'=\sum_{k=1}^j\delta_{y_k}+\omega\wedge\omega'$, where
$\omega\wedge\omega':=\frac12(\omega+\omega' - |\omega-\omega
'|)$. We then have
\begin{eqnarray*}
|F(\omega)-F(\omega')|&\le& |F(\omega)-F(\omega
\wedge\omega')|+|F(\omega')-F(\omega\wedge\omega')|
\\
&\le&\sum_{k=1}^i\Biggl|F\Biggl(\omega\wedge\omega'+\sum
_{l=1}^k\delta_{x_l}\Biggr)-F\Biggl(\omega\wedge\omega'+\sum
_{l=1}^{k-1}\delta_{x_l}\Biggr)\Biggr|
\\
&&{} +\sum_{k=1}^j\Biggl|F\Biggl(\omega\wedge\omega'+\sum
_{l=1}^k\delta_{y_l}\Biggr)-F\Biggl(\omega\wedge\omega'+\sum
_{l=1}^{k-1}\delta_{y_l}\Biggr)\Biggr|
\\
&\leq &\sum_{k=1}^i\varphi(x_k)+\sum_{k=1}^j\varphi(y_k)=\int_E\varphi
\,\mathrm{d}|\omega-\omega'|=d_\varphi(\omega, \omega'),
\end{eqnarray*}
which implies that $\| F\|_{\mathrm{Lip}(d_\varphi)}\leq1$.
\end{pf}

\begin{remarks} When $\varphi=1$, we denote $d_\varphi$ by $d$.
Obviously, $d(\omega, \omega')=|\omega-\omega'|(E)=\|
\omega-\omega'\|_{\mathrm{TV}}$, that is, $d$ is exactly the total variation
distance.
\end{remarks}

The following result, due to the fourth-named author \cite{Wu00},
was obtained by means of the $L^1$-log-Sobolev inequality and
will play an important role.

\begin{lem}[(\cite{Wu00}, Proposition 3.2)]\label{Wuptrf}
Let $F\in L^1(\Omega, P^0)$. If there is some $0\le\varphi\in
L^2(E, m)$ such that $|D_xF(\omega)|\leq\varphi(x)$, $m\otimes
P^0$-a.e., then for any $\lambda\geq0$,
\begin{eqnarray*}
\mathbb{E}^{P^0}
\mathrm{e}^{\lambda(F-P^0(F))}\leq\exp\biggl\{\int_E(\mathrm{e}^{\lambda\varphi
}-\lambda\varphi-1)\,\mathrm{d}m\biggr\}.
\end{eqnarray*}
In particular, if $m$ is finite and $|D_xF(\omega)|\leq1$ for
$m\times P^0$-a.e.~$(x,\omega)$ on $E\times\Omega$ (i.e.,
$\varphi(x)=1)$, then
\begin{eqnarray*}
\mathbb{E}^{P^0}
\mathrm{e}^{\lambda(F-P^0(F))}\leq\exp\{(\mathrm{e}^{\lambda}-\lambda
-1)m(E)\}.
\end{eqnarray*}
\end{lem}

We now state our main result on the Poisson space.
\begin{theorem}\label{main} Let
$(\Omega, \mathcal{F}, P^0)$ be the Poisson space with intensity measure
$m(dx)$ and $\varphi$ a bounded continuous function on $E$ such that
$0<\varphi\leq M$ and $\sigma^2=\int_E \varphi^2 \,\mathrm{d}m<+\infty$. Then
%
\begin{equation}\label{W1honpoisson} \frac1M h_c (W_{1,d_\varphi}(Q,
P^0))\leq
H(Q|P^0)\qquad  \forall Q\in\mathcal{M}_1(\Omega),
\end{equation}
where $c=\sigma^2/M$ and
%
\begin{equation}\label{main2}
h_{c}(r)=c\cdot h\biggl(\frac{r}{c}\biggr), \qquad
h(r)=(1+r)\log(1+r)-r.
\end{equation}
\end{theorem}

Note that $h^*(\lambda):=\sup_{r\ge0}(\lambda r-h(r))=\mathrm{e}^\lambda
-\lambda-1$ and $h_c^*(\lambda)=ch^*(\lambda)$.

\begin{pf*}{Proof of Theorem \ref{main}}
Since the function
$(\mathrm{e}^{\lambda\varphi}-\lambda\varphi-1)/\varphi^2$ is
increasing in $\varphi,$ it is easy to see that
%
\begin{equation}\label{legen}
\int_E(\mathrm{e}^{\lambda\varphi}-\lambda\varphi-1)\,\mathrm{d}m\leq\frac
{\mathrm{e}^{\lambda
M}-\lambda M-1}{M^2}\int\varphi^2\,\mathrm{d}m.
\end{equation}
Further, the Legendre
transformation of the right-hand side of \eqref{legen}
is, for $r\ge0$,
\begin{eqnarray*}
\sup_{\lambda\geq0}\biggl\{\lambda r-\frac{\mathrm{e}^{\lambda M}-\lambda
M-1}{M^2}\int\varphi^2\,\mathrm{d}m\biggr\}&=&\biggl(\frac{r}{M}+\frac{\int
\varphi^2
\,\mathrm{d}m}{M^2}\biggr)\log\biggl(\frac{Mr}{\int\varphi^2 \,\mathrm{d}m}+1
\biggr)-\frac{r}{M}
\\
&=& \frac1M h_{c}(r).
\end{eqnarray*}
The desired result then follows from Theorem \ref{Gozleo}, by Lemma
\ref{Wuptrf}.
\end{pf*}

\begin{remarks} Let
$\beta(\lambda):=\int_E(\mathrm{e}^{\lambda\varphi}-\lambda\varphi-1)\,\mathrm{d}m$ and
$\alpha(r):=\sup_{\lambda\ge0}(\lambda r -\beta(\lambda))$. The
proof above gives us
\[
\alpha(W_{1,d_\varphi}(Q,P^0))\le H(Q|P^0)
\qquad \forall Q\in\mathcal{M}_1(\Omega).
\]
This less explicit inequality is
sharp. Indeed, assume that $E$ is compact and let $F(\omega):=\int_E
\varphi(x) (\omega- m)(\mathrm{d}x)$. We have $\|F\|_{\mathrm{Lip}(d_\varphi)}=1$ and
\[
\log\mathbb{E}^{P^0} \mathrm{e}^{\lambda F} = \beta(\lambda).
\]
The sharpness is then ensured by Theorem \ref{Gozleo}.
\end{remarks}

\begin{prop}\label{prosharp} If $\varphi=1$ and $m$ is finite, then the
inequality \eqref{W1honpoisson} turns out to be
%
\begin{equation}\label{sharp}h_{m(E)}(W_{1, d}(Q, P^0))\le H(Q|P^0)\qquad
\forall
Q\in\mathcal{M}_1(\Omega).
\end{equation}
In particular, for the Poisson measure
$\mathcal{P}(\lambda)$ with parameter $\lambda>0$ on $\mathbb{N}$
equipped with the
Euclidean distance $\rho$,
%
\begin{equation}\label{sharp2}h_\lambda(W_{1,
\rho}(\nu, \mathcal{P}(\lambda)))\le H(\nu|\mathcal{P}(\lambda
)) \qquad \forall\nu
\in
\mathcal{M}_1(\mathbb{N}).
\end{equation}
\end{prop}

\begin{pf} The inequality (\ref{sharp}) is a particular case of
(\ref{W1honpoisson}) with $\varphi=1$ and it holds on $\Omega^0:=\{
\omega\in\Omega; \omega(E)<+\infty\}$
(for $P^0$ is actually
supported in $\Omega^0$ as $m$ is finite). For (\ref{sharp2}), let
$m(E)=\lambda$ and
consider the mapping $\Psi\dvtx \Omega^0\to
\mathbb{N}$, $\Psi(\omega)=\omega(E)$. Since $|\Psi(\omega)-\Psi
(\omega
')|=|\omega(E)-\omega'(E)|\le d(\omega,
\omega')$, $\Psi$ is Lipschitzian with the
Lipschitzian coefficient less than $1$. Thus, (\ref{sharp2}) follows
from (\ref{sharp}) by
\cite{DGW}, Lemma 2.1 and its proof.
\end{pf}

\begin{remarks}\label{rem1} The transportation inequality (\ref{sharp2})
was shown by Liu \cite{Liu09} by means of a tensorization technique
and the approximation of $\mathcal{P}(\lambda)$ by binomial
distributions. It
is optimal (therefore, so is (\ref{sharp})). In fact, consider another
Poisson distribution $\mathcal{P}(\lambda')$ with parameter
$\lambda'>\lambda$. On the one hand,
\begin{eqnarray*} H(\mathcal{P}(\lambda')|\mathcal{P}(\lambda
))&=&\int_\mathbb{N}\log
\frac{\mathrm{d}\mathcal{P}(\lambda')}{\mathrm{d}\mathcal{P}(\lambda)} \,\mathrm{d}\mathcal
{P}(\lambda')=
\sum_{n=0}^{\infty}\mathcal{P}(\lambda')(n)\log\biggl(\frac
{\mathrm{e}^{-\lambda
'}\lambda'^n}{n!}\Big/\frac{\mathrm{e}^{-\lambda}\lambda^n}{n!}\biggr)
\\
&=&\lambda-\lambda'+\sum_{n=0}^{\infty}\mathcal{P}(\lambda')(n)
n\log\frac{\lambda'}{\lambda}
\\
&=&\lambda-\lambda'+\lambda'\log\frac{\lambda'}{\lambda}.
\end{eqnarray*}
On the other hand, let $r:=\lambda'-\lambda>0$. Let $X, Y$ be two
independent random variables having distributions $\mathcal
{P}(\lambda)$ and $\mathcal{P}(r)$, respectively. Obviously, the law
of $X+Y$ is $\mathcal{P}(\lambda').$ Then
\[
W_{1,
\rho}(\mathcal{P}(\lambda'), \mathcal{P}(\lambda))\le\mathbb
{E}|X-(X+Y)|=\mathbb{E}Y=r.
\]
Now,
supposing that $(X, X')$ is a coupling of $\mathcal{P}(\lambda')$ and
$\mathcal{P}(\lambda)$, we have
\[
\mathbb{E}|X-X'| \ge|\mathbb{E}X-\mathbb{E}X'|=r,
\]
which
implies that $W_{1, \rho}(\mathcal{P}(\lambda'), \mathcal
{P}(\lambda))\ge r.$ Then
$W_{1, \rho}(\mathcal{P}(\lambda'), \mathcal{P}(\lambda))= r$ (and
$(X,X+Y)$ is an
optimal coupling for $\mathcal{P}(\lambda)$ and $\mathcal{P}(\lambda
')$). Therefore,
\[
h_{\lambda}(W_{1, \rho}(\mathcal{P}(\lambda'),
\mathcal{P}(\lambda)))=h_{\lambda}(r)=H(\mathcal{P}(\lambda
')|\mathcal{P}(\lambda
)).
\]
Namely, $h_\lambda$ is the optimal deviation function for the Poisson
distribution $\mathcal{P}(\lambda)$.
\end{remarks}

\section{A discrete spin system}\label{sec3}

\textit{The model and the Dobrushin interdependence coefficient.} Let
$\Lambda=\{1, \dots, N\}$ ($2\le N\in\mathbb{N}$) and $\gamma\dvtx
\Lambda\times\Lambda\mapsto[0,+\infty]$ be a \textit{non-negative}
interaction function satisfying $\gamma_{ij}=\gamma_{ji} $ and
$\gamma_{ii}=0$ for all $i,j\in\Lambda$. Consider the Gibbs measure
$P$ on $\mathbb{N}^{\Lambda}$ with
%
\begin{equation}\label{Gibbs}P(x_1, \dots,
x_N)=\mathrm{e}^{-\sum_{i<j}\gamma_{ij}x_ix_j}\prod_{i=1}^{N}\mathcal
{P}(\delta_i)(x_i)\Big/C,
\end{equation}
where $\mathcal{P}(\delta_i)(x_i)=\mathrm{e}^{-\delta_i}
\frac{\delta_i^{x_i}}{x_i!}, x_i\in\mathbb{N}$, is the Poisson distribution
with parameter $\delta_i>0$ and $C$\vspace*{1pt} is the normalization constant.
Here and hereafter, the convention that $0\cdot\infty=0$ is used.
Let $P_i(\mathrm{d}x_i|x_{\Lambda})$ be the given regular conditional
distribution of $x_i$ given $x_{\Lambda\setminus\{i\}},$ which is,
in the present case, the Poisson distribution $\mathcal{P}(\delta_i
\mathrm{e}^{-\sum_{j\neq i}\gamma_{ij}x_j})$ with parameter $\delta_i
\mathrm{e}^{-\sum_{j\neq i}\gamma_{ij}x_j}$, with the convention that the
Poisson measure $\mathcal{P}(0)$ with parameter $\lambda=0$ is the Dirac
measure $\delta_0$ at $0$. Define the Dobrushin interdependence
matrix $C:=(c_{ij})_{i,j\in\Lambda}$ w.r.t.~the Euclidean metric
$\rho$ by
%
\begin{equation}\label{Dobrushindef} c_{ij}=\sup_{x_\Lambda
=x'_\Lambda
{\rm off }j}\frac{W_{1, \rho}(P_i(\mathrm{d}x_i|x_{\Lambda}),
P_i(\mathrm{d}x_i'|x_{\Lambda}') )} {|x_j-x_j'|} \qquad\forall i,
j\in\Lambda
\end{equation}
(obviously, $c_{ii}=0$). The Dobrushin
uniqueness condition \cite{Dobrushin68,Dobrushin70} is then
\[
D:=\sup_{j}\sum_{i}c_{ij}<1.
\]
For this model, we can identify
$c_{ij}.$

\begin{lem}\label{Dobrushin} Recall that $\gamma_{ij}\ge0$. We have
\[
c_{ij}=\delta_i(1-\mathrm{e}^{-\gamma_{ij}}).
\]
\end{lem}

\begin{pf} By Remark \ref{rem1}, if $x_\Lambda=x'_\Lambda$ off $j$,
then
\[
W_{1, \rho}(P_i(\mathrm{d}x_i|x_{\Lambda}),
P_i(\mathrm{d}x_i'|x_{\Lambda}'))=\delta_i
|\mathrm{e}^{-\sum_{k}\gamma_{ik}x_k}-\mathrm{e}^{-\sum_{k}\gamma_{ik}x_k'}|.
\]
Without loss of generality, suppose that $x_j=x_j'+x$ with $x\ge1$.
We have then
\begin{eqnarray*} c_{ij}&=&
\delta_i \sup_{x_\Lambda= x_\Lambda' \mathrm{off} j}\frac
{|\mathrm{e}^{-\sum_{k}\gamma_{ik}x_k}-\mathrm{e}^{-\sum_{k}\gamma
_{ik}x_k'}|}{|x_j-x_j'|}
\\
&=& \delta_i \sup_{x\ge1}\frac{1-\mathrm{e}^{-\gamma_{ij}x}}{x}
\qquad \mbox{(taking $x_k=x_k'=0$ for $k\ne j$, $x_j'=0$)}
\\
&=&\delta_i(1-\mathrm{e}^{-\gamma_{ij}}).
\end{eqnarray*}
Here, the first equality holds since $\gamma_{ij}$ is non-negative and the
last equality is due to the fact that $(1-\mathrm{e}^{-\gamma_{ij}x})/x$ is
decreasing in $x>0.$
\end{pf}

 \textit{The transportation inequality $W_1H$ for
mixed measure.} We return to the general framework of the
\hyperref[sec1]{Introduction}. Let $\mathcal{X}$ be a general Polish space and $d$ be a metric
on $\mathcal{X}$ which is lower semi-continuous on $\mathcal{X}\times
\mathcal{X}$. Consider
a mixed probability measure $\mu:=\int_I \mu_\lambda \,\mathrm{d}\sigma
(\lambda)$ on $\mathcal{X}$, where, for each $\lambda\in I$, $\mu
_{\lambda}$
is a probability on $\mathcal{X}$ and $\sigma$ is a probability
measure on
another Polish space $I$. Let $\rho$ be a lower semi-continuous
metric on $I$.

\begin{prop}\label{translate} Suppose that:
\begin{longlist}[(iii)]
\item[(i)] for any $\lambda\in I$, $\mu_\lambda$ satisfies
$\alpha$--$W_1H$ with deviation function $\alpha\in\mathcal{C}$,
\[
\alpha(W_{1, d}(\nu, \mu_{\lambda}))\le H(\nu|\mu_{\lambda}) \qquad
\forall\nu\in\mathcal{M}_1(\mathcal{X});
\]

\item[(ii)] $\sigma$ satisfies a $\beta$--$W_1H$ inequality on $I$ with
deviation function $\beta\in\mathcal{C}$,
\[
\beta(W_{1,\rho}(\eta,\sigma)) \le H(\eta|\sigma) \qquad \forall
\eta\in\mathcal{M}_1(I);
\]
\item[(iii)] $\lambda\to\mu_\lambda$ is Lipschitzian, that is, for some
constant $M>0$,
\[
W_{1, d}(\mu_{\lambda}, \mu_{\lambda'})\le M \rho(\lambda,
\lambda') \qquad \forall\lambda,\lambda'\in I.
\]
\end{longlist}
The mixed
probability $\mu=\int_{I} \mu_{\lambda} \,\mathrm{d}\sigma(\lambda)$ then satisfies
%
\begin{equation}\label{w1hmunu} \tilde{\alpha}(W_{1, d}(\nu, \mu
))\le
H(\nu|\mu)  \qquad \forall\nu\in\mathcal{M}_1(\mathcal{X}),
\end{equation}
where
\[
\tilde{\alpha}(r)=\sup_{b\ge0}\{b r - [\alpha^*(b)+\beta
^*(b M)]\},\qquad r\ge0.
\]
\end{prop}

\begin{pf} By Gozlan and Leonard's Theorem \ref{Gozleo}, it is
enough to show
that for any Lipschitzian function $f$ on $\mathcal{X}$ with
$\|f\|_{\mathrm{Lip}(d)}\le1$ and $b\ge0$,
\[
\int_\mathcal{X}\mathrm{e}^{b [f(x) - \mu(f)]} \,\mathrm{d}\mu(x)\le\exp\bigl(\alpha
^*(b) +
\beta^*(b M)\bigr).
\]
Let $g(\lambda):=\int_\mathcal{X}f(x) \,\mathrm{d}\mu_\lambda(x)=\mu_\lambda
(f)$. We
have $\sigma(g)=\mu(f)$ and, by Kantorovitch's duality equality and
our condition (iii), $|g(\lambda)-g(\lambda')|\le M
\rho(\lambda,\lambda')$. Using Theorem \ref{Gozleo} and our
conditions (i) and (ii), we then get, for any $b\ge0$,
\begin{eqnarray*} \int_\mathcal{X}\mathrm{e}^{b [f(x) - \mu(f)]} \,\mathrm{d}\mu&=&
\int_I \biggl(\int_\mathcal{X}\mathrm{e}^{b [f(x) - \mu_\lambda(f)]}
\,\mathrm{d}\mu_\lambda(x)\biggr)\mathrm{e}^{b [g(\lambda) - \sigma(g)]}
\,\mathrm{d}\sigma(\lambda),
\\
&\le& \mathrm{e}^{\alpha^*(b)+\beta^*(b M)}
\end{eqnarray*}
the desired result.
\end{pf}

We now turn to a mixed Poisson distribution,
%
\begin{equation}\label{mu}
\mu=\int_{0}^{a}\mathcal{P}(\lambda)\sigma(\mathrm{d}\lambda),
\end{equation}
where $a>0$.
By Proposition \ref{prosharp}, we know that w.r.t.~the Euclidean
metric $\rho$,
\[
h_{\lambda}(W_{1, \rho}(\nu, \mathcal{P}(\lambda)))\le
H(\nu|\mathcal{P}(\lambda))
\]
and $W_{1, \rho}(\mathcal{P}(\lambda), \mathcal{P}(\lambda
'))=|\lambda-\lambda'|.$
Since $h_{\lambda}$ is decreasing in $\lambda,$ the hypotheses in
Proposition \ref{translate} with $E=\mathbb{N}$, $I=[0,a]$, both equipped
with the Euclidean metric $\rho$, are satisfied with
$\alpha(r)=h_a(r)=a h(\frac{r}{a})$ and $\beta(r)=2r^2/a^2$ (the
well-known CKP inequality). On the other hand, obviously,
\[
h(r)=(1+r)\log(1+r)-r\le\frac{r^2}{2}, \qquad r\ge0,
\]
which implies
that
\[
h_{a^2/4}(r)=\frac{a^2}4 h\biggl(\frac{4r}{a^2}\biggr)\le
\frac{2r^2}{a^2}=\beta(r).
\]
Since $h_c^*(\lambda)=c (\mathrm{e}^\lambda-\lambda-1)$,
\[
\sup_{b\ge
0}\{br-[(h_a(b))^{\ast}+(h_{a^2/4}(b))^{\ast}]\}=\sup
_{b\ge
0}\{br-(a+a^2/4)(\mathrm{e}^b-b-1)\}=h_{a+a^2/4}(r).
\]
By
Proposition \ref{translate}, we have, for the mixed Poisson measure
$\mu$ given in (\ref{mu}),
%
\begin{equation}\label{mixed} h_{a+a^2/4}(W_{1,
d}(\nu, \mu))\le H(\nu|\mu)\qquad \forall\nu\in\mathcal
{M}_1(\mathbb{N}
).
\end{equation}

See Chafai and Malrieu \cite{CM09} for fine analysis of
transportation or functional inequalities for mixed measures. We can
now state the main result of this section.
\begin{theorem}\label{dispoisson} Let $P$ be the Gibbs measure given in
\eqref
{Gibbs} with $\gamma_{ij}\ge0$. Assume
Dobrushin's uniqueness condition
\[
D:=\sup_{j\in\Lambda} \sum_{i\in\Lambda}
\delta_i(1-\mathrm{e}^{-\gamma_{ij}})<1.
\]
For any probability measure $Q$ on $\mathbb{N}^\Lambda$ equipped with
the metric $\rho_H(x_\Lambda, y_\Lambda):=\sum_{i\in\Lambda}
|x_i-y_i|$ (the index $H$ refers to Hamming), we then have, for
$c:=\sum_{i\in\Lambda} (\delta_i+\delta_i^2/4)$,
\[
h_c\bigl((1-D)W_{1, \rho_H}(Q, P)\bigr)
\le H(Q|P)\qquad \forall Q\in\mathcal{M}_1(\mathbb{N}^\Lambda).
\]
\end{theorem}

This result, without the extra constants $\delta_i^2/4$, would
become sharp if $\gamma=0$ (i.e., without interaction) or
$P=\mathcal{P}(\delta)^{\otimes\Lambda}$.

\begin{pf*}{Proof of Theorem \ref{dispoisson}} By Theorem
\ref{Gozleo}, it is equivalent to prove that for any
$1$-Lipschitzian functional $F$ w.r.t.~the metric $\rho_H$,
%
\begin{equation}\label{disGL} \log\mathbb{E}^P \mathrm{e}^{\lambda(F-\mathbb
{E}^ P F)}\le
h^*_{c}\biggl(\frac{\lambda}{1-D}\biggr)=c
h^*\biggl(\frac{\lambda}{1-D}\biggr)\qquad \forall\lambda>0.
\end{equation}
We
prove the inequality \eqref{disGL} by the McDiarmid--Rio martingale
method (as in \cite{DGW,Wu06}). Consider the martingale
\[
M_0=\mathbb{E}^{P}(F), \qquad M_k(x_1^k)=\int F(x_1^k, x_{k+1}^N) P
(\mathrm{d}x_{k+1}^N|x_1^k), \qquad 1\le k\le N,
\]
where $x_i^j=(x_k)_{i\le
k\le j}, P(dx_{k+1}^N|x_1^k)$ is the conditional distribution of
$x_{k+1}^N$ given $x_1^k.$ Since $M_N=F,$ we have
\[
\mathbb{E}^P \mathrm{e}^{\lambda(F-\mathbb{E}^P F )}=\mathbb{E}^P\exp
\Biggl(\lambda\sum
_{k=1}^N (M_k-M_{k-1})\Biggr).
\]
By induction, for \eqref{disGL}, it suffices to establish that for
each $k=1, \dots, N, P$-a.s.,
%
\begin{equation}\label{subdis}
\log\int
\exp\bigl(\lambda\bigl(M_k(x_1^{k-1},
x_k)-M_{k-1}(x_1^{k-1})\bigr)\bigr)P(\mathrm{d}x_k|x_1^{k-1})\le
(\delta_k+\delta_k^2/4) h^*\biggl(\frac{\lambda}{1-D}\biggr).
\end{equation}
By
(\ref{mixed}), $P(\mathrm{d}x_k|x_1^{k-1})$, being a convex combination of
Poisson measures $P_k(\mathrm{d}x_k|x_\Lambda)=\mathcal{P}(\delta_k
\mathrm{e}^{-\sum_{j\neq k}\gamma_{kj}x_j})$ (over $x_{k+1}^N$), satisfies
the $W_1H$-inequality with the deviation function
$h_{\delta_k+\delta_k^2/4}$. Hence, by Theorem \ref{Gozleo},
\eqref{subdis} holds if
%
\begin{equation}\label{Lip}|M_k(x_1^{k-1},
x_k)-M_{k}(x_1^{k-1}, y_k)|\le\frac{1}{1-D} |x_k-y_k|.
\end{equation}
In
fact, the inequality \eqref{Lip} has been proven in \cite{Wu06}, step
2 in
the proof of Theorem 4.3. The proof is thus complete.
\end{pf*}

\begin{remarks} For a previous study on transportation inequalities for
Gibbs measures on discrete sites, see Marton \cite{Marton04} and
Wu \cite{Wu06}. Our method here is quite close to that in
\cite{Wu06}, but with two new features: (1) $W_1H$ for mixed
probability measures; (2) Gozlan and L\'eonard's Theorem \ref{Gozleo} as
a new tool.
\end{remarks}

\begin{remarks} Every Poisson distribution $\mathcal{P}(\lambda)$
satisfies the
Poincar\'e inequality (\cite{Wu00}, Remark~1.4)
\[
\operatorname{Var}_{\mathcal{P}(\lambda)} (f) \le\lambda\int_\mathbb{N}(Df(x))^2
\,\mathrm{d}\mathcal{P}(\lambda)(x) \qquad \forall f\in L^2(\mathbb{N},\mathcal
{P}(\lambda)),
\]
where $Df(x):=f(x+1)-f(x)$ and $\operatorname{Var}_\mu(f):=\mu(f^2)-[\mu(f)]^2$ is
the variance of $f$ w.r.t.~$\mu$. By \cite{Wu06}, Theorem 2.2 we
have the following Poincar\'e inequality for the Gibbs measure $P$:
if $D<1$, then
\[
\operatorname{Var}_P(F) \le\frac{\max_{1\le i\le N}\delta_i}{1-D}
\int_{\mathbb{N}^\Lambda} \sum_{i\in\Lambda} (D_i F)^2(x) \,\mathrm{d}P(x) \qquad
\forall
F\in L^2(\mathbb{N}^\Lambda, P),
\]
where $D_iF(x_1,\dots,x_N):=F(x_1,\dots, x_{i-1}, x_i+1,
x_{i+1},\dots, x_N)-F(x_1,\dots,x_N)$. We remind the reader
that an important open question is to prove the $L^1$-log-Sobolev
inequality (or entropy inequality)
\[
H(F P|P) \le C \int_{\mathbb{N}^\Lambda} \sum_{i\in\Lambda} D_i F
\cdot D_i
\log F \,\mathrm{d}P  \qquad \mbox{for all $P$-probability densities } F
\]
(which is equivalent to the exponential convergence in entropy of
the corresponding Glauber system) under Dobrushin's uniqueness
condition, or at least for high temperature.
\end{remarks}

\section{$W_1H$-inequality for the continuum Gibbs measure}\label{sec4}

We now generalize the result for the discrete sites Gibbs measure in
Section~\ref{sec3} to the continuum Gibbs measure (continuous gas model), by an
approximation procedure.

Let $(\Omega,\mathcal{F},P^0)$ be the Poisson space over a compact subset
$E$ of $\mathbb{R}^d$ with intensity $m(\mathrm{d}x)=z \,\mathrm{d}x$, where the Lebesgue
measure $|E|$ of $E$ is positive and finite, and $z>0$ represents
the \textit{activity}. Given a \textit{non-negative} pair-interaction
function $\phi\dvtx \mathbb{R}^d\mapsto[0,+\infty]$, which is
measurable and
even over $\mathbb{R}^d$, the corresponding Poisson space is denoted by
$(\Omega,\mathcal{F}, P^0)$ and the associated Gibbs measure is given by
\[
P^{\phi}(\mathrm{d}\omega)=\frac{\mathrm{e}^{-(1/2)\sum_{x_i, x_j\in{\rm \operatorname{supp}}(\omega), i\ne j}\phi(x_i-x_j)
-\sum_{k, x_i\in{\rm \operatorname{supp}}(\omega)}\phi(x_i-y_k)}}{Z}P^0(\mathrm{d}\omega),
\]
where $Z$ is the
normalization constant and $\{y_k, k\}$ is an at most countable
family of points in $\mathbb{R}^d\backslash E$ such that $\sum_k
\phi(x-y_k)<+\infty$ for all $x\in E$ (boundary condition). The main
result of this section is
the following theorem.

\begin{theorem}\label{w1hgibbs} Assume that the Dobrushin uniqueness condition
holds, that is,
%
\begin{equation}\label{D}D:=z
\int_{\mathbb{R}^d}\bigl(1-\mathrm{e}^{-\phi(y)}\bigr)\,\mathrm{d}y<1.
\end{equation}
Then, w.r.t.~the total
variation distance $d=d_\varphi$ with $\varphi=1$ on $\Omega$,
%
\begin{equation}\label{w1hgibbs2}
h_{z |E|}\bigl((1-D)W_{1,d}(Q,
P^{\phi})\bigr)\le H(Q| P^{\phi}) \qquad \forall Q\in\mathcal
{M}_1(\Omega
).
\end{equation}
\end{theorem}

\begin{remarks} Without interaction (i.e., $\phi=0$), $D=0$ and the
$W_1H$-inequality (\ref{w1hgibbs2}) is exactly the optimal
$W_1H$-inequality for the Poisson measure $P^0$ in Proposition
\ref{prosharp}. In the presence of non-negative interaction $\phi$,
it is well known that $D<1$ is a sharp condition for the
analyticity of the pressure functional $p(z)$: indeed, the radius
$R$ of convergence of the entire series of $p(z)$ at $z=0$ satisfies
$R \int_{\mathbb{R}^d}(1-\mathrm{e}^{-\phi(y)})\,\mathrm{d}y<1$; see \cite{Ru},
Theorem 4.5.3. The corresponding sharp Poincar\'e
inequality for $P^\phi$ was established in  \cite{Wu04}.
\end{remarks}

\begin{pf*}{Proof of Theorem \ref{w1hgibbs}} We shall establish this sharp
$\alpha$--$W_1H$ inequality for $P^\phi$ by approximation.

By part (b$'$) of Theorem \ref{Gozleo}, it is equivalent to show
that for any $F,G\in C_b(\Omega)$ such that $F(\omega)-G(\omega')\le
d(\omega,\omega'),~ \omega,\omega'\in\Omega$, and for any
$\lambda>0$,

\begin{equation}\label{w1hgibbs01} \log\int_\Omega \mathrm{e}^{\lambda
F}\,\mathrm{d}P^\phi\le
\lambda P^\phi(G) + z|E| h^*\biggl(\frac\lambda{1-D}\biggr),
\end{equation}
where $h^*(\lambda)=\mathrm{e}^\lambda-\lambda-1$.

\textit{Step \textup{1}. $\phi$ is continuous and $\{y_k,k\}$ is finite. } We
want to approximate $P^\phi$ by the discrete sites Gibbs measures
given in the previous section. To this end, assume first that
$\phi$ is continuous ($+\infty$ is regarded as the one-point
compactification of $\mathbb{R}^+$) or, equivalently, that $\mathrm{e}^{-\phi
}\dvtx \mathbb{R}
^d\to
[0,1]$ is continuous with the convention that $\mathrm{e}^{-\infty}:=0$.

For each $N\ge2$, let $\{E_1, \dots, E_N\}$ be a measurable
decomposition of $E$ such that, as $N$ goes to infinity, $\max_{1\le
i\le N} \operatorname{Diam}(E_i)\to0$ and $\max_{1\le i\le N}|E_i|\to0$,
where $|E|$ is the Lebesgue measure of $E$ and ${\rm
Diam}(E_i)=\sup_{x,y\in E_i}|x-y|$ is the diameter of $E_i$. Fix
$x_i^0\in E_i$ for each $i$. Consider the probability measure $P_N$
on $\mathbb{N}^\Lambda$ ($\Lambda:=\{1,\dots,N\}$) given by, for all
$(n_1,\dots,n_N)\in\mathbb{N}^\Lambda$,
\begin{eqnarray*}
P_N(n_1,\dots,n_N) &=&(1/Z)\mathrm{e}^{-(1/2)\sum_{i\neq
j}\phi(x_i^0-x_j^0)n_in_j-\sum_{i,k} \phi(x_i^0-y_k)n_i
}\prod_{i=1}^N\mathcal{P}(z|E_i|)(n_i)
\\
&=&(1/Z') \mathrm{e}^{-\sum_{i<
j}\phi(x_i^0-x_j^0)n_in_j}\prod_{i=1}^N\mathcal{P}(\delta_{N,i})(n_i),
\end{eqnarray*}
where $Z,Z'$ are normalization constants and $\delta_{N,i}=z|E_i| \mathrm{e}^{-
\sum_k\phi(x_i^0-y_k)}\le z|E_i|$.
Consider the mapping
$\Phi\dvtx \mathbb{N}^\Lambda\to\Omega$ given by
\[
\Phi(n_1,\dots,n_N)=\sum_{i=1}^N n_i \delta_{x_i^0}.
\]
$\Phi$ is isometric from $(\mathbb{N}^\Lambda,\rho_H)$ to $(\Omega,d)$,
where $d=d_\varphi$ with $\varphi=1$ (given in Section~\ref{sec2}). Finally, let
$P^N$ be the push-forward of $P_N$ by $\Phi$. It is quite direct to
see that $P^N\to P$ weakly.

The Dobrushin constant $D_N$ associated with $P_N$ is given by
\[
D_N =\sup_j\sum_{i}\delta_{N,i}
\bigl(1-\mathrm{e}^{-\phi(x_i^0-x_j^0)}\bigr)\le\sup_j \sum_{i} z|E_i|
\bigl(1-\mathrm{e}^{-\phi(x_i^0-x_j^0)}\bigr).
\]
When $N$ goes to infinity,
\[
\limsup_{N\to\infty}D_N\le\sup_{y\in
\mathbb{R}^d} z\int_E \bigl(1-\mathrm{e}^{-\phi(x-y)}\bigr) \,\mathrm{d}x = z\int_{\mathbb{R}^d}
\bigl(1-\mathrm{e}^{-\phi(x)}\bigr) \,\mathrm{d}x=D.
\]
Therefore, if $D<1$ and $D_N<1$ for all $N$ large enough, then the
$W_1H$-inequality in Theorem~\ref{dispoisson} holds for $P_N$. By
the isometry of the mapping $\Phi$, $P^N$ satisfies the same
$W_1H$-inequality on $\Omega$ w.r.t.~the metric $d$, which gives us,
by Theorem \ref{Gozleo}(b$'$),
\[
\log\mathbb{E}^{P^N} \mathrm{e}^{\lambda F} \le\lambda P^N(G) +
\biggl(\sum_{i\in\Lambda} [\delta_{N,i} + \delta_{N,i}^2/4]\biggr)
h^*\biggl(\frac\lambda{1-D_N}\biggr).
\]
By letting $N$ go to infinity, this yields (\ref{w1hgibbs01}),
for $P^N\to P^\phi$ weakly and
\[
\sum_{i\in\Lambda} [\delta_{N,i} + \delta_{N,i}^2/4]\le
\sum_{i\in\Lambda} z|E_i|(1+ z|E_i|/4)\to z|E|.
\]

\textit{Step \textup{2}. General $\phi$ and $\{y_k,k\}$ is finite.} For general
measurable non-negative and even interaction function $\phi$, we take
a sequence of continuous, even and non-negative functions $(\phi_n)$
such that $1-\mathrm{e}^{-\phi_n}\to1-\mathrm{e}^{-\phi}$ in $L^1(\mathbb{R}^d,\mathrm{d}x)$. Now,
note that $ \frac{\mathrm{d}P^{\phi_n}}{\mathrm{d}P^0} \to\frac{\mathrm{d}P^{\phi}}{\mathrm{d}P^0} $
in $L^1(\Omega,P^0)$, that is, $P^{\phi_n}\to P^\phi$ in total
variation. Hence, (\ref{w1hgibbs01}) for $P^{\phi_n}$ (proved in step
1) yields (\ref{w1hgibbs01}) for $P^{\phi}$.

\textit{Step \textup{3}. General case.} Finally, if the set of points $\{y_k,k\}$
is infinite, approximating $\sum_{k=1}^\infty\phi(x_i-y_k)$ by
$\sum_{k=1}^n \phi(x_i-y_k)$ in the definition of $P^\phi$, we get
(\ref{w1hgibbs01}) for $P^{\phi}$, as in step 2.
\end{pf*}

\begin{remarks}\label{rem41}
The explicit Poissonian concentration
inequality (\ref{a1}) follows from Theorem \ref{w1hgibbs} by Theorem
\ref{Gozleo}(c) (with $n=1$) by noting that the observable
$F(\omega)=\omega(f)/(2N)^d$ there is Lipschitzian w.r.t.~$d$ with
$\|F\|_{\rm Lip(d)}\le M/(2N)^d$ and $h(r)\ge(r/2)\log(1+r)$.
\end{remarks}

\begin{remarks} A quite curious phenomena occurs in the continuous gas
model: the {\it extra} constant $\delta_i^2/4$ coming from the
mixture of measures now disappears.
\end{remarks}

\section*{Acknowledgments}
We are grateful to the referee for his conscientious comments. Yutao
Ma was supported by NSFC Grant No.~10721091.

\printhistory


\begin{thebibliography}{99}

\bibitem{BG99}
Bobkov, S.G. and G\"{o}tze, F. (1999).
Exponential integrability and transportation cost related to
logarithmic {S}obolev inequalities.
\textit{J. Funct. Anal.} \textbf{163} 1--28.
\MR{1682772}

\bibitem{BV05}
 Bolley, F. and Villani, C. (2005).
Weighted {C}sisz\'ar--{K}ullback--{P}insker inequalities and
applications to transportation inequalities.
\textit{Ann. Fac. Sci. Toulouse} \textbf{14} 331--352.
\MR{2172583}


\bibitem{CM09} Chafai, D. and Malrieu, F. (2010). On fine properties of
mixtures with respect to concentration and Sobolev type
inequalities. \textit{Ann. Inst. H.
Poincar\'e.} Preprint. To appear.

\bibitem{DGW} Djellout, H., Gullin, A. and Wu, L.M. (2004).  Transportation
cost-information inequalities for
random dynamical systems and diffusions. \textit{Ann. Probab.} \textbf{32} 2702--2732.
\MR{2078555}

\bibitem{Dobrushin68} Dobrushin, R.L. (1968).  The description of a random
field by means of
conditional probabilities and condition of its regularity. {\it
Theory Probab. Appl.} \textbf{13} 197--224.
\MR{0231434}

\bibitem{Dobrushin70} Dobrushin, R.L. (1970). Prescribing a system of random
variables by conditional distributions. \textit{Theory Probab. Appl.}
 \textbf{15} 458--486.

\bibitem{GL} Gozlan, N. and L\'eonard, C.  (2007). A large deviation approach
to some transportation cost inequalities. \textit{Probab. Theory
Related Fields} \textbf{139} 235--283.
\MR{2322697}

\bibitem{Kel84} Kellerer, H.G. (1984). Duality theorems for marginal
problems. \textit{Z. Wahrsch. Verw. Gebiete}
\textbf{67} 399--432.
\MR{0761565}

\bibitem{Leo} L\'eonard, C. (2007). Transport inequalities: A large
deviation point of
view. In \textit{Course in Chinese Summer School for Ph.D. Students,
Wuhan}.


\bibitem{Liu09} Liu, W. Optimal transportation-entropy inequalities
for several usual distributions on
$\mathbb{R}$. {Preprint. Submitted.}

\bibitem{Mar96}
 Marton, K. (1996).
Bounding $\bar{d}$-distance by informational divergence: A way to
prove measure concentration.
\textit{Ann. Probab.} \textbf{24} 857--866.
\MR{1404531}

\bibitem{Marton04} Marton, K. (2004). Measure concentration for Euclidean
distance in the case of dependent random variables.
\textit{Ann. Probab.} \textbf{32}  2526--2544.
\MR{2078549}

\bibitem{Ru} Ruelle, D.  (1969). \textit{Statistical Mechanics:
Rigorous Results}. New York:  Benjamin.
\MR{0289084}

\bibitem{Villani} C. Villani. (2003).  \textit{Topics in Optimal Transportation}.
 Providence, {RI}: {Amer. Math. Soc.}
\MR{1964483}

\bibitem{Wu00} Wu, L.M. (2000). A new modified logarithmic Sobolev inequality for
Poisson processes and several applications. \textit{Probab. Theory
Related Fields} \textbf{118} 427--438.
\MR{1800540}

\bibitem{Wu04} Wu, L.M. (2004). Estimate of the spectral gap for continuous
gas. \textit{Ann. Inst. H. Poincar\'e Probab. Statist.} \textbf{40} 387--409.
\MR{2070332}

\bibitem{Wu06} Wu, L.M.  (2006). Poincar\'e and transportation inequalities
for Gibbs
measures under the Dobrushin uniqueness condition. \textit{Ann. Probab.} \textbf{34} 1960--1989.
\MR{2271488}

\end{thebibliography}
\end{document}